\documentclass[12pt, leqno]{article}
\usepackage{amssymb}
\usepackage{amscd}
\usepackage{amsmath}
\usepackage{amsthm}
\usepackage{graphics}
\usepackage{mathptmx}
\usepackage{newtxtext}       %
\usepackage[varvw]{newtxmath}       
\usepackage{mathrsfs}
\usepackage{refcount}
\def\overset#1#2{{\mathrel{\mathop {{#2}_{}}\limits^{#1}}}}
\def\underset#1#2{{\mathrel{\mathop {{}_{} {#2}}\limits_{{#1}_{}}}}}
\def\upplim_#1{\underset{#1}{\overline\lim}\;}
\def\lowlim_#1{\underset{#1}{\underline\lim}\;}
\newtheorem{thm}{Theorem}[section]

\newtheorem{defn}[thm]{\indent{\it Definition}\rm }

\newtheorem{lem}[thm]{Lemma}

\newtheorem{rmk}[thm]{\indent \rm {\it Remark}}

\newcommand{\C}{{\mathbf{C}}}
\newcommand{\sC}{{\mathscr{C}}}

\newcommand{\fa}{\forall}

\newcommand{\del}{{\partial}}
\newcommand{\delbar}{\bar{\partial}}

\newcommand{\iso}{\cong}

\newcommand{\II}{{I \hskip-1pt I}}

\newcommand{\lto}{\longrightarrow}
\newcommand{\sM}{\mathscr{M}}

\newcommand{\N}{\mathbf{N}}

\renewcommand{\O}{{\mathscr{O}}}

\newcommand{\vp}{\varphi}

\newcommand{\refn}{\getrefnumber}

\newcommand{\supp}{\mathrm{Supp}\,}

\newcommand{\Z}{\mathbf{Z}}
\newcommand{\wbar}{\widebar}
\newenvironment{pf}{\par{\it Proof. }}{\qed \medskip}
\numberwithin{equation}{section}

\title{
Inverse of Abelian  Integrals  and  
 Ramified Riemann Domains (II): on a Proof of the  Behnke--Stein
Theorem
}
\date{\small The University of Tokyo}
\author{
Junjiro Noguchi\thanks{
  Research supported in part by Grant-in-Aid
 for Scientific Research (C) 25K07025. \hfill\break
}
}
\begin{document}
\setlength{\baselineskip}{16pt}
\maketitle
\begin{abstract}
In the first titled paper (Math.\ Ann.\ 367 (2017), 229--249)
a new proof of the Behnke--Stein theorem was given
using some scalar function which is an analogue of Hartogs' radius.
In the proof the vanishing of $H^2(X, \Z)$ for an open Riemann surface
$X$ was used. It is a most common way to refer some references on
algebraic topology for that vanishing.
The aim of the present (II) is to provide a short analytic proof
of the vanishing of $H^2(X, \Z)$
using the {\it Continuous Cousin Problem} on  $X$,
which is self-contained from the viewpoint of complex analysis.
The proof is independent from algebraic topology.
\end{abstract}

\thispagestyle{empty}

\section{Introduction and main results}

It is a well-known fundamental theorem of
Behnke--Stein \cite{BS} that every open Riemann surface $X$ is Stein.
Actually they proved for $X$ the Runge approximation theorem,
Mittag-Leffler's theorem, and Weierstrass' theorem, which imply
the Stein conditions; the converse holds by the general
Oka--Cartan theory (cf. e.g., \cite{nog1}).

In the first titled paper \cite{nog2} \S1.2.2
we gave a new proof of the Behnke--Stein theorem by showing the
Stein conditions for $X$.
In the proof we used some scalar function which is an analogue of
 Hartogs' radius. After knowing the Steinness of a bounded
domain $\Omega \Subset X$, the scalar function is used to show
that a pair $(\Omega_1, \Omega_2)$ of domains
 $\Omega_1 \Subset \Omega_2 (\Subset X)$ is Runge if there is no
relatively compact connected component of $\Omega_2 \setminus \Omega_1$
in $\Omega_2$.
 To define the scalar function we used the triviality
of the canonical bundle $K_X \iso \mathbf{1}_X$,
 for which  a topological fact of the vanishing of the \v{C}ech
cohomology $H^2(X, \Z)$ was referred:
 It is  common  to cite some books on
algebraic topology for that vanishing, but in fact the proof of it
is not easy, rigorously rather involved using the duality of
 the cohomologies with
 the homologies and singular homologies; cf., e.g., 
Ahlfors--Sario \cite{AS} Chap.\ 1 for triangulations,
homologies and singular homologies,
Voisin \cite{vo} Chap.\ 4 \S4.3, Theorem 4.47
for the isomorphism between \v{C}ech and singular cohomologies,
and Hatcher \cite{ha} Chap.\ 3 \S3.1, Theorem 3.2
for the duality by the universal coefficient theorem.

The aim of the present (II) is to provide a short analytic
proof of the vanishing of $H^2(X, \Z)$
 using the {\it Continuous Cousin Problem} 
 introduced in \cite{nog3} \S3.4.2, \cite{nog4} \S8.4.1,
which is self-contained from the viewpoint of complex analysis.

\begin{lem}\label{h20}
If $X$ is a Stein Riemann surface, then $H^2(X, \Z)=0$.
\end{lem}

This is an intermediate result, by which
combined with \cite{nog2} Theorem 1.10 we complete the proof of:
\begin{thm}\label{bst}
Every open Riemann surface $X$ is Stein and satisfies $H^2(X, \Z)=0$.
\end{thm}

In the present method, it is a point to claim the Steinness and
the vanishing of $H^2(X,\Z)$ at the same time.
The proof is independent from algebraic topology. 

{\it Acknowledgment.} The author is very grateful to Shuya Matsumoto for
 his interest in the present subject and for providing
 references on algebraic topology.

\section{Preparations}
We recall the notion of ``{\it Continuous Cousin}.
In general, let $X$ be a connected complex manifold of dimension $n$
 (or a complex space as well) with structure sheaf $\O_X$.
Then $X$ is called a {\it Stein manifold} if the following conditions
are satisfied (cf.\ Grauert--Remmert \cite{GR} Introduction \S4):
\begin{enumerate}\setlength{\itemsep}{-2pt}
\item
(Second countability) $X$ satisfies the second countability condition.
\item
(Holomorphic convexity)  For every relatively compact subset $E \Subset X$
the holomorphic hull $\hat{E}_X:=\{p \in X: |f(p)| \leq 
\sup_E |f|, \, \fa f \in \O(X)\}$ is compact, where
$\O(X)$ denotes the space of all holomorphic functions on $X$.
\item
(Holomorphic separation) For any distinct point $p,q \in X$
there is a holomorphic function $f \in \O(X)$ such that $f(p) \not= f(q)$.
\item
(Global holomorphic local coordinates) For every point
$p \in M$ there are holomorphic functions
 $f_j \in \O(X), 1 \leq j \leq n$ such that $(f_1, \ldots, f_n)$
is a holomorphic local coordinate system about $p$.
\end{enumerate}

Let $\{U_\alpha\}_{\alpha \in \Gamma}$ be an open covering of $X$,
and let $\vp_\alpha \in C^0(U_\alpha)$ be a continuous function on
each $U_\alpha$.
\begin{defn}\label{ccd}\rm
\begin{enumerate}\setlength{\itemsep}{-2pt}
\item
The family $\{(U_\alpha, \vp_\alpha)\}_{\alpha\in\Gamma}$
of pairs $(U_\alpha, \vp_\alpha)$ is called {\it continuous Cousin data}
on $X$ if $\vp_\alpha - \vp_\beta \in \O(U_\alpha \cap U_\beta)$
 for all $\alpha, \beta \in \Gamma$.
\item
A continuous function $\vp \in C^0(X)$ is called a 
{\it solution} of the continuous Cousin data
$\sC_c:=\{(U_\alpha, \vp_\alpha)\}_{\alpha\in\Gamma}$
  if $\vp-\vp_\alpha$ is holomorphic
in every $U_\alpha$; if it exists, $\sC_c$ is said to be solvable.
\end{enumerate}
\end{defn}

\begin{thm}[\textbf{Continuous Cousin}]\label{cct}
If $X$ is Stein, then every continuous Cousin data on $X$ is solvable.
\end{thm}
\begin{pf}
The same arguments as in the proof of \cite{nog2} Theorem 3.4.22
in the domain case in $\C^n$ work on  Stein $X$.
\end{pf}

 In fact, the proof is similar to that of Cousin I Problem with
replacing ``meromorphic'' with ``continuous'', but the above
solves simultaneously the Cousin I and II Problems, and 
the $\delbar$-equation with aid of Dolbeault's lemma
on Stein manifolds (cf.\ \cite{nog3} Chap.\ 3, \cite{nog4} Chap.\ 8).

Let $\sC_{\II}=\{(U_\alpha, \psi_\alpha)\}_{\alpha \in \Gamma}$
be Cousin II data, i.e., $\psi_\alpha \in \sM^*(U_\alpha)$ is
a meromorphic function in each $U_\alpha$ such that
$\psi_\alpha$ is not constantly zero in any connected component
of $U_\alpha$ and
$\psi_\alpha/\psi_\beta \in \O^*(U_\alpha \cap U_\beta)$ is
a zero-free holomorphic function in $U_\alpha \cap U_\beta$
for every pair of $\alpha$ and $\beta$.
We say that $\sC_{\II}$ is {\it analytically solvable}
(resp.\ {\it topologically solvable}) if there is
a meromorphic function (resp.\ continuous function) $\psi$ on $X$
(resp.\ a dense open set of $X$)
such that $\psi/\psi_\alpha \in \O^*(U_\alpha \cap U_\beta)$
(resp.\ $\psi/\psi_\alpha \in C^*(U_\alpha \cap U_\beta)$)
for all $\alpha, \beta \in \Gamma$, where $C^*(U_\alpha \cap U_\beta)$
denotes the set of all zero-free continuous functions in 
$U_\alpha \cap U_\beta$.
We call $\psi$ an {\it analytic} (resp.\ {\it topological})
solution of $\sC_{\II}$.

Then the Oka Principle is stated as follows (cf.\ \cite{nog3} Theorem
3.4.51):
\begin{thm}[\textbf{Oka Principle}]\label{okp}
Assume that $X$ is a Stein manifold.
Then Cousin II data  on $X$ is analytically solvable
if it is topologically solvable.
\end{thm}
\begin{pf}
For a convenience we give a proof using the continuous Cousin Theorem \refn{cct}.
Let $\sC_{\II}=\{(U_\alpha, \psi_\alpha)\}_{\alpha \in \Gamma}$
be Cousin II data on $X$. Assume that there is a topological solution
 $\psi$ of $\sC_{\II}$.
We may assume that every $U_\alpha$ is simply connected.
Since $\psi_\alpha/\psi \in C^*(U_\alpha)$,
there is a $1$-valued branch $f_\alpha :=\log (\psi_\alpha/\psi)$
in $U_\alpha$. Then $\sC_I:=\{(U_\alpha, f_\alpha)\}$
is continuous Cousin data.
By Theorem \refn{cct} there is a solution $f \in C^0(X)$
of $\sC_I$, so that $g_\alpha:=f- \log {\psi_\alpha/\psi}
\in \O(U_\alpha)$. Then we have
\[
 \psi_\alpha e^{g_\alpha}=\psi e^f.
\]
The expression of the right-hand side is global and
the left-hand side is a meromorphic function locally
with the given zeros and poles.
Therefore $\vp:=\psi_\alpha e^{g_\alpha}$ in each $U_\alpha$
is an analytic solution of $\sC_\II$ on $X$. 
\end{pf}

Let $L \to X$ be a holomorphic line bundle over $X$, and
denote by $\O(L)$ the sheaf of germs
of holomorphic sections in $L$ over $X$.
Then $\O(L)$ is a coherent sheaf over $X$.
We write $\O(L)(U)$ for the space of all holomorphic sections in $L$
over an open set $U$ of $X$.
By the Oka theory or the Oka--Cartan theory we have:
\begin{thm}\label{hq0}
If $X$ is Stein, then $H^q(X, \O(L))=0$   for $q \geq 1$.
\end{thm}

For the proof, cf., e.g., the arguments in \cite{nog2} Theorem 3.4.22,
or \cite{nog1} Theorem 4.4.2.

\section{Proof of Lemma \refn{h20}}
(1) As usual, we write $\O^*_X$ for the sheaf of multiplicative
groups of units of $\O^*_X$.
It follows from the short exact sequence
\[
 0 \to \Z \lto \O_X ~ \overset{e^{2\pi i *}}{\lto} ~ \O^*_X \to 0,
\]
and Theorem \refn{hq0} that
\[
 0=H^1(X, \O_X) \lto H^1(X, \O^*_X) ~ \overset{\iso}{\lto}~ 
H^2(X, \Z) \to H^2(X, \O_X)=0.
\]
Therefore it suffices to show that $H^1(X, \O^*_X)=0$. i.e.,
 every holomorphic line bundle $L \to X$ is trivial over $X$.
We are going to prove $H^1(X, \O^*_X)=0$
without assuming Weierstrass' theorem; here, rather oppositely,
we would like to derive it from Theorem \refn{okp}.

Let $L \to X$ be a holomorphic line bundle over $X$.
Let $\sM(L)$ denote the sheaf of germs of meromorphic
sections in $L$ over $X$ and let $\sM(L)(U)$ denote
the space of sections in $\sM(L)$ over an open set $U$ of $X$.
The following lemma reduce the problem to the divisor case.
\begin{lem}\label{hsec}
There is a meromorphic section 
$\phi \in \sM(L)(X)$  such that $\phi \not\equiv 0$.
\end{lem}
\begin{pf}
Let $p_0 \in X$ be any fixed point and let
 $z$ be a local coordinate in a neighborhood
$U_0$ about $p_0$ such that $z(p_0)=0$ and
 $L|_{U_0} \iso U_0 \times \C$ is trivial.
Set $g_0=1/z$ to be a meromorphic section in $L$
on $U_0$ and let $g_1=0$ to be a zero section in $L$
on $U_1:=X \setminus \{p\}$.
Then $g_{01}=g_0|_{U_0 \cap U_1} \in \O(L)(U_0 \cap U_1)$, and
 $(\{U_j\}_{j=0.1}, \{g_{01})\}$ is a $1$-cocycle
with coefficient in $\O(L)$.
It follows from Theorem \refn{hq0} that there is a
meromorphic section $\phi \in \sM(L)(X)$ over $X$ with pole of order $1$ at $p_0$,
so that $\phi \not\equiv 0$.
\end{pf}

Set $m_0=-1$ and let $p_\nu \in X ~(1 \leq \nu \leq N \leq \infty)$ be
zeros of $\phi$ above with order $m_\nu \in \N$.
We set a divisor on $X$:
\begin{equation}\label{div}
\delta=\sum_{\nu=0}^N m_\nu p_\nu .
\end{equation}

We are going to obtain a meromorphic function $F \in \sM^*(X)$
whose divisor is $\delta$. If it exists, then $\frac1F \phi$ is
a zero-free holomorphic section in $L$, and hence $L$ is trivial.

The rest of the arguments is very elementary and we proceed step by step.

(2)   Let $\Delta_0=\{|z|< r\} ~(r>1)$ be a local coordinate disk in $X$,
 and let $a_j \in \{|z|<1\}, j=0,1$, be two distinct points.
Then we have the  point divisors $a_j, j=0,1$,  on $X$.
In general, we denote by $\sC_\II(\vartheta)$ the Cousin II data
associated with a divisor $\vartheta$ on $X$; in the present case,
we consider the Cousin II data $\sC_\II(a_0 - a_1)$
 consisting of those data such that
$z-a_0$ about $a_0$, $\frac1{z-a_1}$ about $a_1$, and $1$
in the other points in $X$. Here ``$a_0 - a_1''$ does not
represents an operation in complex numbers, but in divisors on $X$.

\begin{lem}\label{point1}
Let the notation be as above. Then  there is a topological
solution $g_{(a_0.a_1)}$ of  $\sC_\II(a_0-a_1)$
such that $g_{(a_0, a_1)}|_{X \setminus \{|z| \leq r'\}} \equiv 1$
 with $1< r' <r$.
\end{lem}
\begin{pf}
It is noticed that a branch of $\log\frac{z-a_0}{z-a_1}$ is determined
as a $1$-valued holomorphic
function in a neighborhood of the circle $|z|=1$.
Let $\lambda(x)$ be a continuous function on $X$ such that
\begin{enumerate}\setlength{\itemsep}{-2pt}
\item
  $0 \leq \lambda(x) \leq 1$,
\item
 $\lambda(z)=1$ in $|z|<1$,
\item
 $\supp \lambda \subset \{|z|<r \}$.
\end{enumerate}
We set
\begin{equation}
g_{(a_0, a_1)}= \begin{cases}\vspace{5pt}
\frac{z-a_0}{z-a_1} & \quad \hbox{in } |z|<1, \\
e^{\lambda(z)\log \frac{z-a_0}{z-a_1}} & \quad \hbox{in }
\{ 1 \leq |z|< r\}, \\
1 & \quad \hbox{in } X \setminus \supp \lambda
\end{cases}
\end{equation}
Then $g_{(a_0,a_1)}$ is a topological solution of
$\sC_\II(a_0-a_1)$.
There is an $r'<r$ with $\supp \lambda \subset \{|z|< r'\}$.
Then $g_{(a_0, a_1)}|_{X \setminus \{|z| \leq r'\}} \equiv 1$.
\end{pf}

Let $b \in X$ be a point different to $a_0$ and let $C$
be a Jordan curve from $a_0$ to $b$.
\begin{lem}\label{point2}
 Let $V \supset C$ be an open subset of $X$.
Then there is a topological solution $g_{(a_0, b)}$ of the Cousin II
data $\sC_\II(a_0-b)$ such that
$g_{(a_0, b)}|_{X \setminus V} \equiv1$.
\end{lem}
\begin{pf}
We take a sequence of points on $C$,
\[
 a_0, ~ a_1, ~ \ldots, ~ a_l=b
\]
such that $a_{j -1}$ and $a_j$ are contained in a same
local coordinate disk $\Delta_j \Subset V$ for
every $j$ with $1 \leq j \leq l$.
It follows from Lemma \ref{point1} that for every such $j$
there is a topological solution $g_{(a_{j-1}, a_j)}$
of the Cousin II data 
$\sC_\II (a_{j-1} - a_j))$
such that $g_{(a_{j-1}, a_j)}|_{X\setminus \Delta_j}\equiv 1$.
Set
\[
 g_{(a_0, b)}=\prod_{j=1}^{l} g_{(a_{j-1}, a_j)}.
\]
Then the intermediate zeros and singularities of $g_{(a_{j-1}, a_j)}$
are canceled out,  so that $g_{(a_0, b)} \in C^*(X\setminus\{a_0, b\})$
and it is a topological solution of $\sC_\II(a_0 -b)$.
Since $\bigcup_{j=1}^{l} \Delta_j \Subset V$,
we have that $g_{(a_0, b)}|_{X \setminus V}\equiv1$.
\end{pf}

(3) Now we go back to \eqref{div}.
First we take an increasing connected open covering of $X$,
\[
 \Omega_1 \Subset \Omega_2 \Subset \cdots,
\quad \bigcup_{\mu=1}^\infty \Omega_\mu=X,
\]
such that for every $\mu \geq 1$,
\begin{enumerate}\setlength{\itemsep}{-2pt}
\item
there is  no relatively compact connected component in
 $X \setminus \wbar{\Omega}_\mu$,
\item
$\{p_\nu: \nu \geq 0\} \cap \del \Omega_\mu=\emptyset$,
where $\del \Omega_\mu$ is the boundary of $\Omega_\mu$.
\end{enumerate}

We may assume that
 $\{p_\nu: \nu \geq 0 \}\cap \Omega_1=\{p_\nu: 0 \leq \nu \leq \nu_1\}$
with $\nu_1 \in \N$.
For every $p_\nu ~(0 \leq \nu \leq \nu_1)$
we take $q_{\nu 2} \in \Omega_2 \setminus \wbar\Omega_1$.
Then there is a Jordan curve from $p_\nu$ to $q_{\nu2}$ in $\Omega_2$.
It follows from Lemma \ref{point2} that for every such $\nu$
 there is a topological solution $g_{(p_\nu, q_{\nu2})}$ of
 $\sC_\II(p_\nu-q_{\nu2})$
with $g_{(p_\nu, q_{\nu2})}|_{X\setminus \Omega_2}\equiv 1$.
Set
\begin{equation}\label{ts1}
G_1 = g_{(p_0 , q_{02})}^{m_0}\cdot g_{(p_1, q_{12})}^{m_1} \cdot \cdots
\cdot g_{(p_{\nu_1}, q_{{\nu_1}2})}^{m_{\nu_1}}.
\end{equation}
Then $G_1$ is a topological solution of
$\sC_\II \big(\sum_{\nu=0}^{\nu_1} m_\nu a_\nu -
\sum_{\nu=0}^{\nu_1}m_\nu q_{\nu 2} \big)$
such that $G_1|_{X \setminus \Omega_2} \equiv 1$.

Let
\[
 \{p_\nu: \in \N\}\cap (\Omega_2 \setminus \wbar\Omega_1)=
\{p_{\nu_1+1}, p_{\nu_1+2}, \ldots, p_{\nu_2}\}.
\]
For each $p_\nu ~(\nu_1+1 \leq \nu \leq \nu_2)$ we take
a point $q_{\nu 3} \in \Omega_3 \setminus \wbar\Omega_2$ in the connected component
of $\Omega_3 \setminus \wbar\Omega_1$, containing $p_\nu$.
Then there is a Jordan curve from $p_\nu$ to $q_{\nu 3}$ in
 $\Omega_3 \setminus \wbar\Omega_1$.
By Lemma \ref{point2} there is a topological solution 
$g_{(p_\nu, q_{\nu3})}$ of $\sC_\II(p_\nu -q_{\nu 3})$
such that 
$g_{(p_\nu, q_{\nu3})}|_{X \setminus (\Omega_3 \setminus\Omega_1)} \equiv1$;
in particular, $g_{(p_\nu, q_{\nu3})}|_{\Omega_1} \equiv1$.

Furthermore, for every $q_{\nu 2}$ with $0 \leq \nu \leq \nu_1$
we take a point $q_{\nu 3} \in \Omega_3 \setminus \wbar\Omega_2$
in the connected component
of $\Omega_3 \setminus \wbar\Omega_1$, containing $q_{\nu2}$.
By Lemma \ref{point2} there is a topological solution 
$g_{(q_{\nu2}, q_{\nu3})}$ of $\sC_\II(q_{\nu 2} -q_{\nu 3})$
such that 
$g_{(q_{\nu2}, q_{\nu3})}|_{X \setminus (\Omega_3 \setminus\Omega_1)} \equiv1$;
in particular.
$g_{(q_{\nu2}, q_{\nu3})}|_{\Omega_1} \equiv1$.
We set
\begin{equation}\label{ts2}
G_2=G_1 \times \prod_{\nu=0}^{\nu_1} g_{(q_{\nu 2}, q_{\nu 3})}^{m_\nu}
\times \prod_{\nu=\nu_1+1}^{\nu_2} g_{(p_\nu, q_{\nu 3})}^{m_\nu}.
\end{equation}
Then $G_2$ is a topological solution of 
$\sC_\II \big(\sum_{\nu=0}^{\nu_2} m_\nu a_\nu -
\sum_{\nu=0}^{\nu_2}m_\nu q_{\nu 3} \big)$
such that
\begin{align*}
 G_2|_{\Omega_1} &= G_1|_{\Omega_1}, \\
 G_2|_{X \setminus \Omega_3} &\equiv 1.
\end{align*}

We proceed by induction to construct a topological solution
$G_\nu$ of $\sC_\II \big(\sum_{\nu=0}^{\nu_\mu} m_\nu a_\nu -
\sum_{\nu=0}^{\nu_\mu} m_\nu q_{\nu \mu+1} \big)$
such that
\begin{align*}
 G_\mu|_{\Omega_{\mu-1}} &= G_{\mu-1}|_{\Omega_{\mu-1}}, \\
 G_\mu|_{X \setminus \Omega_{\mu+1}} &\equiv 1.
\end{align*}

Finally we set
\begin{equation}\label{tsf}
G(x)= \lim_{\mu \to \infty} G_\mu(x), \quad x \in X.
\end{equation}
Then $G|_{\Omega_\mu}=G_\mu|_{\Omega_\mu}$ for all $\mu \in \N$
and $G$ is a topological solution of
$\sC_\II (\delta)$.
 We have proved:
\begin{lem}
There is a topological solution of $\sC_\II(\delta)$.
\end{lem}

It follows from Theorem \refn{okp} (Oka Principle) that
there is a meromorphic function $F \in \sM^*(X)$, whose divisor is $\delta$.

Thus we have completed the proof of Lemma \ref{h20}.
\qed

\begin{rmk}
The above proof to obtain a meromorphic function $F$ works
for any divisor $\delta$ on $X$.
\end{rmk}

\rightline{Graduate School of Mathematical Sciences}
\rightline{University of Tokyo}
\rightline{Komaba, Meguro-ku, Tokyo 153-8914}
\rightline{Japan}
\rightline{e-mail: noguchi@ms.u-tokyo.ac.jp}
\end{document}